\documentclass{amsart}

\usepackage{amssymb,amsfonts, latexsym,amsmath,hhline,array,longtable}
\usepackage{pdfsync,color, comment,colortbl, mathrsfs,stmaryrd,cite,graphicx,yfonts}

\input xy
\xyoption{all}

\usepackage{ifpdf}
\ifpdf \usepackage[colorlinks=true, citecolor=blue, linkcolor=blue, urlcolor=blue]{hyperref} \fi

\newcommand{\cal}{\mathcal}

\newtheorem{formula}{}[section]
\newtheorem{definition}[formula]{Definition}
\newtheorem{corollary}[formula]{Corollary}
\newtheorem{remark}[formula]{Remark}
\newtheorem{lemma}[formula]{Lemma}
\newtheorem{theorem}[formula]{Theorem}
\newtheorem{proposition}[formula]{Proposition}
\newtheorem*{claim}{Claim}

\def\thrm{\begin{theorem*}}
\def\thrml#1{\begin{theorem}\label{#1}}
\def\ethrm{\end{theorem}}
\def\prpstn{\begin{proposition}}
\def\prpstnl#1{\begin{proposition}\label{#1}}
\def\eprpstn{\end{proposition}}
\def\rmrk{\begin{remark}}
\def\rmrkl#1{\begin{remark}\label{#1}}
\def\ermrk{\end{remark}}
\def\dfntn{\begin{definition}}
\def\dfntnl#1{\begin{definition}\label{#1}}
\def\edfntn{\end{definition}}
\def\nmrt{\begin{enumerate}}
\def\enmrt{\end{enumerate}}
\def\tm#1{\item[{\rm (#1)}]}
\def\qtnl#1{\begin{equation}\label{#1}}
\def\eqtn{\end{equation}}
\def\lmm{\begin{lemma}}
\def\lmml#1{\begin{lemma}\label{#1}}
\def\elmm{\end{lemma}}
\def\crllr{\begin{corollary}}
\def\crllrl#1{\begin{corollary}\label{#1}}
\def\ecrllr{\end{corollary}}
\def\css{\begin{cases}}
\def\ecss{\end{cases}}
\def\prf{\begin{proof}}
\def\eprf{\end{proof}}
\def\clm{\begin{claim}}
\def\eclm{\end{claim}}

\def\cH{{\cal H}}

\def\cX{{\cal X}}
\def\cY{{\cal Y}}

\def\fK{{\mathfrak K}}
\def\fM{{\mathfrak M}}

\def\fS{{\mathfrak S}}
\def\fs{{\mathfrak s}}

\def\fT{{\mathfrak T}}
\def\ft{{\mathfrak t}}
\def\fX{{\mathfrak X}}

\def\sC{{\mathscr{C}}}

\def\uni{{\textgoth 1}}

\DeclareMathOperator{\aut}{Aut}

\DeclareMathOperator{\id}{id}

\DeclareMathOperator{\inv}{Inv}

\DeclareMathOperator{\iso}{Iso}

\DeclareMathOperator{\mon}{Mon}

\DeclareMathOperator{\pr}{pr}

\DeclareMathOperator{\secc}{Sec}

\DeclareMathOperator{\rad}{rad}

\DeclareMathOperator{\res}{res}
\DeclareMathOperator{\rk}{rk}

\DeclareMathOperator{\sym}{Sym}
\DeclareMathOperator{\WL}{WL}

\def\bone{{\bf 1}}
\def\grp#1{\langle {#1}\rangle}

\def\phmaa#1{{\phantom{x}\hspace{-2mm}^{#1}}}
\def\qaq{\quad\text{and}\quad}

\def\wh{\widehat}

\def\VRTB#1{*=<10mm>[o][F-]{#1}}

\begin{document}

\title{On the Weisfeiler-Leman dimension of circulant graphs}
\author{Yulai Wu$^*$}
\address{School of Mathematics and Statistics, Hainan University, Haikou, China}
\email{wuyl@hainanu.edu.cn}
\thanks{$^*$Partly supported by Hainan Province Natural Science Foundation of China, grant No. 120RC452.}
\author{Ilia Ponomarenko$^\dag$}
\address{School of Mathematics and Statistics, Hainan University, Haikou, China and Steklov Institute of Mathematics at St. Petersburg, Russia}
\email{inp@pdmi.ras.ru}
\thanks{$^\dag$Supported by Natural Science  Foundation of China, grant No. 12361003.}

\begin{abstract}
A circulant graph is a Cayley graph of a finite cyclic group. The Weisfeiler-Leman-dimension of a circulant graph $X$ with respect to the class of all circulant graphs is the smallest positive integer~$m$ such that the $m$-dimensional Weisfeiler-Leman algorithm correctly tests the isomorphism between $X$ and any other circulant graph. 	It is proved that for  a circulant graph of order $n$  this dimension is less than or equal to $\Omega(n)+3$, where $\Omega(n)$ is the total number of prime divisors of~$n$. 
\end{abstract}

\maketitle

\section{Introduction}

The Weisfeiler-Leman algorithm was initially proposed in \cite{WLe68} for testing graph isomorphism. It is in this capacity that a multidimensional variant of this algorithm was used in  Babai's fundamental result \cite{Babai2019} on the graph isomorphism problem. For a positive integer $m$ and an arbitrary graph $X$ (directed or undirected), the $m$-dimensional Weisfeiler-Leman algorithm (the $m$-dim $\WL$ algorithm) canonically constructs a special coloring of the $m$-tuples consisting of vertices of~$X$. Two given graphs are declared isomorphic if the resulting colorings do not distinguish these graphs; in this case the graphs are said to be \emph{$\WL_m$-equivalent}. The algorithm works correctly if the $\WL_m$-equivalence of the input graphs implies that they are isomorphic. In many cases, the $m$-dim $\WL$ algorithm works correctly already for small $m$, for example, if one of the input graphs is planar, one can take $m=3$ \cite{Kiefer2017}. However, there are pairs of $n$-vertex graphs with arbitrary large~$n$, for which $m$ cannot be taken less than $cn$ for some constant $0<c< 1$, see ~\cite{CaiFI1992}.

A natural formalism closely related with the above discussion was proposed in monograph~\cite{Grohe2017}. A  \emph{Weisfeiler-Leman dimension} (\emph{$\WL$-dimension}) of a graph $X$ is defined to be the smallest positive integer~$m$ such that any graph $X'$ which is $\WL_m$-equivalent to the graph~$X$ is isomorphic to~$X$. More or less recent results about the $\WL$-dimension of graphs in some special classes can be found in, e.g.,~\cite{Grohe2021}. When $X$  belongs to a fixed class $\fK$ of graphs and $X'$ runs over~$\fK$, we speak of the $\WL$-dimension of $X$ with respect to~$\fK$. 
 
In the present paper we focus on the following question: how large can be the $\WL$-dimension of a circulant graph, i.e., a Cayley graph of a cyclic group. It is known that the $\WL$-dimension of a random circulant graph is at most~$2$ \cite{Verbitsky2024}, and that of every circulant graph of prime power order  is at most~$3$ \cite{Ponomarenko2022a}; the latter upper bound is tight. At present, it is not clear whether it is possible to bound the $\WL$-dimension of an arbitrary circulant graph by a constant not depending on the order of a graph. Our main result (Theorem~\ref{260723a} below) shows that the $\WL$-dimension of any  circulant graph of order~$n$ with respect to the class of all circulant graphs is at most~$O(\log n)$.

\thrml{260723a}
The $\WL$-dimension of any  circulant graph of order~$n$ with respect to the class of all circulant graphs is at most~$\Omega(n)+3$, where $\Omega(n)$ is the total number of prime divisors of $n$. 
\ethrm

It seems that the bound $\Omega(n)+3$ in Theorem~\ref{260723a} is not tight. At least, it can be reduced to $3$ if $n$ is a power of a prime. Moreover, it is not known whether the $\WL$-dimension of any  circulant graph  with respect to the class of all circulant graphs can be bounded from above by a constant. Finally, it is not clear whether Theorem~\ref{260723a} remains true if one removes the phrase ``with respect to the class of all circulant graphs''.

Let us discuss the main ideas in the proof of Theorem~\ref{260723a}. First, we note that each graph $X$ is associated with a coherent configuration $\cX=\WL(X)$ which can be thought as an arc-colored graph of a special form (the exact definition of coherent configuration is given in Section~\ref{070324a}). The concept of the $\WL$-dimension is naturally extended to coherent configurations (see Section~\ref{070324b}), and then we need to estimate the $\WL$-dimension  $\dim_{\scriptscriptstyle\WL}(\cX)$ of the coherent configuration~$\cX$. In this case,  $\cX$ is  a Cayley scheme over a cyclic group or, in other terminology, a circulant scheme, and we can take advantage of the theory of such schemes, developed in series of papers \cite{EvdP2003,EvdP2004a,EvdP2012a,Evdokimov2016,Evdokimov2015}. The relevant details of this theory are presented in Section~\ref{070324c} and Subsection~\ref{280124a}.

When passing from graphs to coherent configurations, the condition of $\WL_m$-equivalence of two graphs $X$ and $X'$ transforms into the condition of the existence of an algebraic isomorphism $\varphi:\WL(X)\to \WL(X')$, the restriction of which to any $m$-vertex subset is induced by  isomorphism (Lemma~\ref{310723a}). By Muzychuk’s theorem (Theorem~\ref{160823c}), algebraically isomorphic circulant schemes are isomorphic and therefore we may assume that $\WL(X)=\WL(X')$, and that $\varphi$ is an algebraic automorphism of the scheme $\cX=\WL(X)$.

In a special class of quasinormal circulant schemes, each algebraic isomorphism is uniquely defined by a multiplier (Lemma~\ref{160823a}), that is, a system of consistent automorphisms of some sections of the underlying cyclic group. The set of these sections is uniquely determined by the scheme $\cX$ and the cardinality of this set does not exceed $\Omega(n)$. A careful analysis of the multiplier corresponding to the algebraic automorphism $\varphi$ shows that if $m\ge \Omega(n)+3$, then $\varphi$ is induced by the isomorphism (Theorem~\ref{170823h}). Thus if the scheme $\cX$ is quasinormal, then $\dim_{\scriptscriptstyle\WL}(X)=\dim_{\scriptscriptstyle\WL}(\cX)\le \Omega(n)+3$.

A reduction from arbitrary circulant schemes to quasinormal ones is established in Section~\ref{100324a}. The key point here is  that any non-quasinormal scheme $\cX$ can be extended to a larger circulant scheme $\cX^{\star}$ the structure of which is, in a sense, under control. From Theorem~\ref{220124b}, it follows  that $\dim_{\scriptscriptstyle\WL}(\cX)\le \dim_{\scriptscriptstyle\WL}(\cX^{\star})$. Since the scheme $\cX^{\star}$ is strictly larger than $\cX$, we arrive at quasinormal case after finitely many extensions.

From the above discussion, it should be clear that the proof of our main result heavily uses coherent configurations, some tools to pass from $\WL_m$-equivalent graphs to algebraic isomorphisms, and  circulant schemes. To make the paper as selfcontained as possible, we cite  relevant facts in Sections~\ref{070324a}, \ref{070324b}, and \ref{070324c}.

\section{Coherent configurations}\label{070324a}
In the further presentation of the basics of  coherent configurations, we follow the monograph~\cite{CP2019}; all facts not explained in detail can be found there.

\subsection{Notation}
Throughout the paper, $\Omega$ denotes a finite set. For a set $\Delta\subseteq \Omega$, the Cartesian product $\Delta\times\Delta$ and its diagonal are denoted by~$\bone_\Delta$ and $1_\Delta$, respectively. If $\Delta=\{\alpha\}$, we abbreviate $1_\alpha:=1_{\{\alpha\}}$. For a relation $s\subseteq\bone_\Omega$, we set $s^*=\{(\alpha,\beta): (\beta,\alpha)\in s\}$, $\alpha s=\{\beta\in\Omega:\ (\alpha,\beta)\in s\}$ for all $\alpha\in\Omega$, and define $\grp{s}$ as the minimal (with respect to inclusion) equivalence relation on~$\Omega$, containing~$s$. The \emph{support} of~$s$ is defined to be the smallest set $\Delta\subseteq\Omega$ such that $s\subseteq\bone_\Omega$.

For any collection $S$ of relations, we denote by $S^\cup$ the set of all unions of elements of~$S$, and consider $S^\cup$ as a poset with respect to inclusion. The \emph{dot product} of relations $r,s\subseteq\Omega\times\Omega$ is denoted by $$r\cdot s=\{(\alpha,\beta)\!:\ (\alpha,\gamma)\in r,\ (\gamma,\beta)\in s\text{ for some }\gamma\in\Omega\}.$$ 

The set of classes of  an equivalence relation $e$ on~$\Omega$ is denoted by $\Omega/e$. For $\Delta\subseteq\Omega$, we set $\Delta/e=\Delta/e_\Delta$, where $e_\Delta=\bone_\Delta\cap e$. If the classes of $e_\Delta$ are singletons, $\Delta/e$ is identified with $\Delta$. Given a relation $s\subseteq \bone_\Omega$, we put
\qtnl{130622a}
s_{\Delta/e}=\{(\Gamma,\Gamma')\in \bone_{\Delta/e}: s^{}_{\Gamma,\Gamma'}\ne\varnothing\},
\eqtn
where $s^{}_{\Gamma,\Gamma'}=s\cap (\Gamma\times \Gamma')$, and  abbreviate $s^{}_\Gamma:=s^{}_{\Gamma,\Gamma}$. For a fixed $s$, the set of all  equivalence relations $e$ on~$\Omega$, such that $e\cdot s=s\cdot e=s$, contains  the largest (with respect to inclusion) element which is denoted by  $\rad(s)$ and called the {\it radical} of~$s$. Obviously, $\rad(s)\subseteq\grp{s}$.

For a set $B$ of bijections $f:\Omega\to\Omega'$, subsets $\Delta\subseteq \Omega$ and $\Delta'\subseteq \Omega'$, equivalence relations $e$ and $e'$ on $\Omega$ and  $\Omega'$, respectively, we put
$$
B^{\Delta/e,\Delta'/e'}=\{f^{\Delta/e}:\ f\in B,\ \Delta^f=\Delta',\ e^f=e'\},
$$
where $f^{\Delta/e}$ is the bijection from $\Delta/e$ onto $\Delta'/e'$ induced by~$f$; we also  abbreviate $B^{\Delta/e}:=B^{\Delta/e,\Delta'/e}$ if $\Delta'$ is clear from context.

Under a tuple on $\Omega$, we mean an $m$-tuple $x=(x_1,\ldots,x_m)$ for some natural number~$m$. We put  $\Omega(x):=\{x_1,\ldots,x_m\}$. For a mapping  $f$ from $\Omega$ to another set, we put $x^f= ({x_1}^f,\ldots,{x_m}^f)$ and extend this notation to sets of $m$-tuples in a natural way.

\subsection{Coherent configurations}\label{270224b}
Let $S$ be a partition of $\bone_\Omega$. A pair $\mathcal{X}=(\Omega,S)$ is called a \emph{coherent configuration} on $\Omega$ if 
\nmrt
\tm{CC1}  $1_\Omega\in S^\cup$,
\tm{CC2} $s^*\in S$ for all $s\in S$, 
\tm{CC3} given $r,s,t\in S$, the number $c_{rs}^t=|\alpha r\cap \beta s^{*}|$ does not depend on $(\alpha,\beta)\in t$. 
\enmrt
The number $|\Omega|$ is called the {\it degree} of~$\cX$. A coherent configuration~$\cX$ is said to be {\it trivial} if  $S=S(\cX)$ consists of~$1_\Omega$ and its complement (unless $\Omega$ is a singleton), \emph{discrete} if~$S$ consists of  singletons, {\it homogeneous} or a {\it scheme} if $1_\Omega\in S$.  

With any permutation group $G\le\sym(\Omega)$, one can associate a coherent configuration $\inv(G)=(\Omega,S)$, where $S$ consists of the relations $(\alpha,\beta)^G=\{(\alpha^g,\beta^g): \ g\in G\}$ for all $\alpha,\beta\in\Omega$. A coherent configuration $\cX$ is said to be \emph{regular} if $\cX=\inv(G)$, where the group $G$ is regular, or equivalently, if $|\alpha s|=1$ for all $\alpha\in\Omega$ and $s\in S$.

\subsection{Isomorphisms and algebraic isomorphisms}\label{110622w} A  {\it combinatorial isomorphism} or, briefly, \emph{isomorphism} from a coherent configuration $\cX=(\Omega, S)$ to a coherent configuration $\cX'=(\Omega', S')$ is defined to be a bijection $f: \Omega\rightarrow \Omega'$ such that  $S^f=S'$, where $S^f=\{s^f:\ s\in S\}$. In this case,  $\cX$ and $\cX'$ are said to be {\it isomorphic}; the set of all isomorphisms from $\cX$ to~$\cX'$ is denoted by $\iso(\cX,\cX')$, and by $\iso(\cX)$ if $\cX=\cX'$. The group $\iso(\cX)$ contains a normal subgroup
$$
\aut(\cX)=\{f\in\sym(\Omega):\ s^f=s \text{ for all } s\in S\}
$$
called the {\it automorphism group} of $\cX$.

A bijection $\varphi:S\rightarrow S'$ is called an \emph{algebraic isomorphism} from $\cX$ to $\cX'$ if for all $r,s,t\in S$, we have
$$
c_{\varphi (r)\varphi (s)}^{\varphi(t)}=c_{rs}^t.
$$
In this case, $|\Omega'|=|\Omega|$ and  $1^{}_{\Omega'}=\varphi(1_\Omega)$. The set of all algebraic isomorphisms from $\cX$ to $\cX'$ is denoted by $\iso_{alg}(\cX,\cX')$, and by $\aut_{alg}(\cX)$ if $\cX=\cX'$.

Any isomorphism $f\in\iso(\cX,\cX')$ induces the algebraic isomorphism $\varphi_f:\cX\to \cX'$, $s\mapsto s^f$.  We put
$$
\iso(\cX,\cX',\varphi)=\{f\in\iso(\cX,\cX'):\ \varphi_f=\varphi\}
$$
and abbreviate $\iso(\cX,\varphi):=\iso(\cX,\cX,\varphi)$. Note that $\aut(\cX)=\iso(\cX,\id_\cX)$, where $\id_\cX$ is the trivial (identical) algebraic automorphism of~$\cX$. If the scheme $\cX$ is regular, then $\iso(\cX,\cX',\varphi)\ne\varnothing$ for any~$\varphi$ (see~\cite[Theorem~2.3.33]{CP2019}).

\subsection{Tensor product}
Given two arbitrary coherent configurations $\cX_1=(\Omega_1,S_1)$ and $\cX_2=(\Omega_2,S_2)$, denote by $S_1\otimes S_2$ the set of all relations 
$$
s_1\otimes s_2=\left\{((\alpha_1,\alpha_2),(\beta_1,\beta_2))\in (\Omega_1\times\Omega_2)^2:\ (\alpha_1,\beta_1)\in s_1,\ (\alpha_2,\beta_2)\in s_2\right\},
$$
where $s_1\in S_1$ and  $s_2\in S_2$. Then the pair $\cX_1\otimes\cX_2=(\Omega_1\times \Omega_2,S_1\otimes S_2)$ is a coherent configuration. It is called the
\emph{tensor product} of~$\cX_1$ and~$\cX_2$.

\subsection{Extensions} There is a natural partial ordering\, $\le$\, on the set of all coherent configurations  on~$\Omega$. Namely, given two such coherent configurations~ $\cX$ and $\cY$, we set
$$
\cX\le\cY\quad \Leftrightarrow\quad S(\cX)^\cup\subseteq S(\cY)^\cup,
$$
and say that $\cY$ is the \emph{extension} (or fission) of~$\cX$. The minimal and maximal elements with respect to this order are the \emph{trivial} and {\it discrete} coherent configurations, respectively. It holds true that $\inv(\aut(\cX))\ge \cX$.

The coherent configuration $\WL(X)$ of a graph $X$, mentioned in the Introduction, is
precisely the minimal coherent configuration on the vertex set of $X$ that contains the arc
set of $X$ as the union of basis relations.

\subsection{Relations and fibers}  Let $\cX=(\Omega,S)$ be a coherent configuration. The elements of $S$ and of $S^\cup$ are called {\it basis relations} and \emph{relations} of ~$\cX$, respectively. The basis relation containing a pair $(\alpha,\beta)\in\bone_\Omega$ is denoted by $r_\cX(\alpha,\beta)$ or simply by~$r(\alpha,\beta)$. The set of all relations of $\cX$ is closed with respect to intersections, unions, and the dot product.

Any set $\Delta\subseteq\Omega$ such that $1_\Delta\in S$ is called a {\it fiber} of $\cX$. In view of condition~(CC1), the set~$F=F(\cX)$  of all fibers forms a partition of~$\Omega$. Every basis relation is contained in the Cartesian product of two uniquely determined fibers. 

Every algebraic isomorphism $\varphi\in\iso_{alg}(\cX,\cX')$ can be extended in a natural way to a poset isomorphism $S^\cup\to (S')^\cup$, denoted also by~$\varphi$. This poset isomorphism respects the set theoretical operations and the dot product. Furthermore, $\varphi$ induces a fiber preserving poset isomorphism $F^\cup\to (F')^\cup$, $\Delta\mapsto\Delta^\varphi$, where the set $\Delta^\varphi$ is defined by the equality $\varphi(1_\Delta)=1_{\Delta^\varphi}$. 

Let $\cX\le \cY$ and $\cX'\le\cY'$. We say that the algebraic isomorphism $\varphi$ \emph{is extended} to an algebraic isomorphism $\psi\in\iso_{alg}(\cY,\cY')$ or that $\psi$ \emph{extends} $\varphi$ if $\psi(s)=\varphi(s)$ for all $s\in S$. 

\subsection{Parabolics and quotients}\label{140322w}  A relation $e$ of $\cX$  that is an equivalence relation on the support of~$e$ is called a \emph{partial parabolic} of~$\cX$, and a \emph{parabolic} of~$\cX$ if the support coincides with~$\Omega$. The relations $\grp{s}$ and $\rad(s)$ are partial parabolics for all $s\in S^\cup$. Given a partial parabolic $e$, we put $S_{\Omega/e}=\{s_{\Omega/e}:\ s\in S,\ s_{\Omega/e}\ne\varnothing\}$, and put $S_\Delta=\{s_\Delta:\ s\in S,\ s_\Delta\ne\varnothing\}$ for any $\Delta\in\Omega/e$. Then 
$$
\cX_{\Omega/e}=(\Omega/e,S_{\Omega/e})\qaq \cX_\Delta=(\Delta,S_\Delta)
$$  
are coherent configurations, called a \emph{quotient} of~$\cX$ modulo~$e$ and  a \emph{restriction} of~$\cX$ to~$\Delta$, respectively. We have 
$$
F(\cX_{\Omega/e})=\{\Gamma_{\Omega/e}:\ \Gamma\in F,\ \Gamma_{\Omega/e}\ne\varnothing\}\qaq  F(\cX_\Delta)=\{\Gamma\cap\Delta:\ \Gamma\in F,\ \Gamma\cap\Delta\ne\varnothing\},
$$ 
respectively. In particular, $\cX_{\Omega/e}$ and $\cX_\Delta$ are schemes if so is~$\cX$. 

Let $\varphi\in\iso_{alg}(\cX,\cX')$ be an algebraic isomorphism and  $e$ a (partial) parabolic of~$\cX$. Then  $\varphi(e)$ is a (partial) parabolic of~$\cX'$ with the same number of classes. Moreover, 
\qtnl{030322a}
\varphi(\grp{s})=\grp{\varphi(s)}\qaq \varphi(\rad(s))=\rad(\varphi(s)).
\eqtn
The algebraic isomorphism $\varphi$ induces the algebraic isomorphism
\qtnl{060524a}
\varphi^{}_{\Omega^{}/e^{}}:\cX^{}_{\Omega^{}/e^{}}\to \cX'_{\Omega'/e'},\ s^{}_{\Omega^{}/e^{}}\mapsto s'_{\Omega'/e'},
\eqtn
where $e'=\varphi(e)$ and $s'=\varphi(s)$.  Now let  $\Delta\in\Omega/e$ and $\Delta'\in\Omega/e'$. Assume that there is $\Gamma\in F$ such that $\Delta\cap\Gamma\ne\varnothing\ne \Delta'\cap\Gamma'$, where $\Gamma'=\Gamma^\varphi$. Then  $\varphi$ induces an algebraic isomorphism 
\qtnl{200322u}
\varphi^{}_{\Delta,\Delta'}:\cX^{}_{\Delta^{}}\to \cX'_{\Delta'},\ s^{}_{\Delta^{}}\mapsto s'_{\Delta'}.
\eqtn
This  algebraic isomorphism always exists for all $\Delta$ and $\Delta'$ if $\cX$ (and hence $\cX'$) is a scheme.

\subsection{Sections} Let $\cX$ be a coherent configuration, $e$ a partial parabolic of $\cX$, and $\Delta$ a class of a partial parabolic containing~$e$. The quotient set $\fS=\Delta/e$ is called a {\it section} of~$\cX$. Any element of $\fS$ is of the form $\alpha_\fS:=\alpha e$ for some $\alpha\in\Delta$; we extend this notation by setting  $\alpha_\fS=\varnothing$  for all $\alpha\not\in \Delta$. For any $s\in S$, we define the relation~$s_\fS$ on $\fS$ by formula~\eqref{130622a}.  

The set of all sections of~$\cX$  is denoted by~$\secc(\cX)$. It   is partially ordered, namely, $\Delta/e\preceq \Delta'/e'$,  whenever $\Delta\subseteq\Delta'$ and $e'\subseteq e$. For a section $\fS=\Delta/e\in\secc(\cX)$, we put $\cX_\fS=(\cX_\Delta)_{\Delta/e}$. 

Let $\varphi\in\iso_{alg}(\cX,\cX')$,  $e'=\varphi(e)$, and $\Delta'$  a class of the $\varphi$-image of the partial parabolic of~$\cX$, containing the class~$\Delta$.  The algebraic isomorphisms~\eqref{060524a} and~\eqref{200322u} induce the  algebraic isomorphism 
$$
\varphi_{\fS,\fS'}=(\varphi_{\Delta,\Delta'})_{\Delta/e,\Delta'/e'}
$$
from the coherent configurations  $\cX^{}_\fS$ to the coherent configuration $\cX'_{\fS'}$. 

\subsection{Point extensions}
Let $m\ge 1$ be an integer and $x\in\Omega^m$.
The  {\it point extension}~$\cX_x$ of  the coherent configuration~$\cX$ with respect to  $x$  is defined to be the smallest coherent configuration that is an extension of~$\cX$, and contains the singleton relations $1_{x_1},\ldots,1_{x_m}$.  Clearly, the point extension ~$\cX_x$ depends on the set $\Omega(x)$ only. When the tuple is irrelevant, we use the term ``point extension'',  and ``one-point extension'' if $m=1$.

Let $\varphi\in\iso_{alg}(\cX,\cX')$, and let $x\in\Omega^m$, $x'\in{\Omega'}^m$. An algebraic isomorphism $\psi\in\iso_{alg}(\cX^{}_{x^{}},\cX'_{x'})$  is  called an {\it $(x,x')$-extension } of $\varphi$ if it extends $\varphi$ and
$$
\psi(1_{x^{}_i})=1_{x'_i},\quad i=1,\ldots,m.
$$
Note that the $(x,x')$-extension is unique if it exists. The following statement is straightforward.

\lmml{030923h} 
Let $\cX$ and $\cX'$ be schemes, $\fS\in\secc(\cX)$ and  $\fS'\in\secc(\cX')$, and
$$
\varphi\in\iso_{alg}(\cX,\cX')\qaq\varphi_{\fS,\fS'}\in\iso_{alg}(\cX^{}_{\fS^{}},\cX'_{\fS'}).
$$
If $\varphi$  has the $(x,x')$-extension for some tuples $x$ and~$x'$, then $\varphi_{\fS,\fS'}$ has the $(x^{}_{\fS^{}},x'_{\fS'})$-extension, where $x^{}_{\fS^{}}$ (respectively, $x'_{\fS'}$) is a tuple with nonempty entries $(x_i)_{\fS^{}}$ (respectively, $(x_i')_{\fS'}$).
\elmm

\section{The WL-dimension of coherent configurations}\label{070324b}
In this section, we give the definition of the $\WL$-dimension of coherent configuration in terms of $m$-ary coherent configurations, $m\ge 2$. The details can be found in~\cite{Chen2023a,Chen2023,Ponomarenko2022a}.

\subsection{$m$-ary coherent configurations}
Let  $m\ge 1$ be an integer. For a tuple $x\in\Omega^m$, denote by $\rho(x)$ the equivalence relation on $M=\{1,\ldots,m\}$ such that $(i,j)\in\rho(x)$ if and only if  $x_i=x_j$.   The monoid of all maps from $M$ to itself is denoted by $\mon(M)$. Given $\sigma\in\mon(M)$, we set 
$$
x^\sigma=(x_{1^\sigma},\ldots,x_{m^\sigma}).
$$ 

Let $\fX$ be a partition of~$\Omega^m$. The class of $\fX$, containing $x\in\Omega^m$, is denoted by~$[x]$. Given $X_1,\ldots,X_m\in \fX$, we denote by $n(x;X_1,\ldots,X_m)$ the number of all $\alpha\in\Omega$ such that $x_{i\leftarrow \alpha} \in X_i$ for all~$i\in M$, where
$
x_{i\leftarrow \alpha}=(x_1,\ldots,x_{i-1},\alpha,x_{i+1},\ldots,x_m).
$

A partition $\fX$ of  $\Omega^m$ is called an \emph{$m$-ary coherent configuration} on~$\Omega$  if the following conditions are satisfied for all $X\in \fX$:
\nmrt
\tm{CC1'} $\rho(x)$ does not depend on $x\in X$,
\tm{CC2'} $X^\sigma\in \fX$  for all $\sigma\in\mon(M)$,
\tm{CC3'} for any $X_0,X_1,\ldots,X_m\in \fX$, the number $n_{X_1,\ldots,X_m}^{X_0}=n(x_0;X_1,\ldots,X_m)$ 
does not depend on $x_0\in X_0$.
\enmrt
The $2$-ary coherent configurations are the coherent configurations in the sense of Subsection~\ref{270224b}, and the definitions below are compatible with what introduced there.

Given $X\subseteq \Omega^m$ and  $k\in M$, put $\pr_k X=\{\pr_k(x):\ x\in X\}$, where $\pr_k(x)=(x_1,\ldots,x_k)$. In accordance with \cite{AndresHelfgott2017} (see also \cite{Ponomarenko2022a}), if $\fX$ is an $m$-ary coherent configuration, then 
\qtnl{270123a}
\pr_k \fX=\{\pr_k X:\ X\in\fX\}
\eqtn
is a   $k$-ary coherent configuration on~$\Omega$; it is called the \emph{$k$-projection} of~$\fX$.

Let $\fX'$ be an $m$-ary coherent configuration on~$\Omega'$. An \emph{isomorphism} from~$\fX$ to~$\fX'$ is  a bijection $f: \Omega^m\rightarrow {\Omega'}\phmaa{m}$ such that  $X^f\in \fX'$ and $(X^\sigma)^f=(X^f)^\sigma$ for all $X\in \fX$ and all $\sigma\in\mon(M)$. An \emph{algebraic isomorphism} from $\fX$ to $\fX'$ is a bijection $\varphi:\fX\to\fX'$ such that
\qtnl{110522a}
\varphi(X^\sigma)=\varphi(X)^\sigma\qaq n_{X_1,\ldots,X_m}^{X_0}=n_{\varphi(X_1),\ldots,\varphi(X_m)}^{\varphi(X_0)}
\eqtn
for all $X,X_0,\ldots,X_m\in \fX$ and $\sigma\in\mon(M)$. 
Again, for $m=2$, this definition is compatible with that given in Subsection~\ref{110622w}. For any $k\in M$, the algebraic isomorphism~$\varphi$ induces the  algebraic isomorphism between the $k$-projections,
\qtnl{150223a}
\varphi_k:\pr_k \fX\to\pr_k\fX',\ \pr_k X\mapsto \pr_k \varphi(X).
\eqtn

\subsection{The Weisfeiler-Leman algorithm}
Let $m\ge 2$. Given a coherent configuration $\cX=(\Omega,S)$, the $m$-dim $\WL$ algorithm constructs an $m$-ary coherent configuration $\WL_m(\cX)$ on~$\Omega$, such that 
\qtnl{080924a}
\pr_2 \WL_m(\cX)\ge \cX.
\eqtn
 We say that $\cX$ is  \emph{$\WL_m$-equivalent} to a coherent configuration~$\cX'$ with respect to an algebraic isomorphism $\varphi\in\iso_{alg}(\cX,\cX')$ if there exists an algebraic isomorphism  
$$
\hat\varphi:\WL_m(\cX)\to\WL_m(\cX')
$$ 
such that $\hat\varphi_2(s)=\varphi(s)$ for all $s\in S$ (see~\eqref{080924a}), where $\hat\varphi_2$ is the $2$-projection of~$\hat\varphi$. The $m$-dim $\WL$ algorithm \emph{identifies}~$\cX$ if  for any coherent configuration  $\cX'$, each algebraic isomorphism $\varphi\in\iso_{alg}(\cX,\cX')$ with respect to which  $\cX$ and $\cX'$ are $\WL_m$-equivalent, is induced by an isomorphism. The \emph{$\WL$-dimension} $\dim_{\scriptscriptstyle\WL}(\cX)$ of the coherent configuration~$\cX$ is defined to be the smallest positive integer $m\ge 2$ such that the $m$-dim $\WL$ algorithm identifies~$\cX$.

\subsection{The Cai-F\"urer-Immerman theorem}\label{060223e}
To deal with the $\WL_m$-equivalence, we  recall (very briefly) the definition of the  pebbling game $\sC_m(\varphi)$ for two coherent configurations $\cX$ and $\cX'$ and an algebraic isomorphism $\varphi\in\iso_{alg}(\cX,\cX')$. 

There are two players, called Spoiler and Duplicator, and $m$ pairwise distinct pebbles, each given in duplicate.   The game consists of rounds and each round consists of two parts. At the first part, Spoiler chooses a set~$A'$ of points in~$\cX$ or in~$\cX'$. Duplicator responds with a set $A$ in the other coherent configuration, such that $|A|=|A'|$ (if this is impossible, then Duplicator loses). At the second part, Spoiler places one of the pebbles\footnote{It is allowed to move previously placed pebbles to other 	vertices and place more than one pebble on the same vertex.}  on a point in~$A$. Duplicator responds  by placing the copy of the pebble on some point of~$A'$. 

The configuration  after a round is determined by a bijection $f:\Delta\to\Delta'$, where $\Delta\subseteq\Omega$ (respectively, $\Delta'\subseteq\Omega' $) is the set of points  in $\cX$ (respectively, in~$\cX'$), covered by pebbles, and any two points $\alpha\in\Delta$ and $\alpha^f\in\Delta'$ are covered by the copies of the same pebble. Duplicator wins the round if 
$$
(s_{\Delta^{}})^f=\varphi(s)_{\Delta'}
$$
for all $s\in S$ such that $s_\Delta\ne\varnothing$. Spoiler wins if Duplicator  does not win.

The game starts from an initial configuration (which  is considered as the configuration after zero rounds), i.e.,  a pair $(x,x')\in\Omega^k\times {\Omega'}\phmaa{k}$  for which $k\le m$ and the points~$x^{}_i$ and~$x'_i$ are covered with copies of the same pebble, $i=1,\ldots,k$. We say that Duplicator  has a winning strategy  for the game $\sC_m(\varphi)$  on $\cX$ and $\cX'$ with  initial configuration $(x,x')$ if, regardless of Spoiler's actions, Duplicator wins after any number of rounds.

\lmml{120324a}
Let $T^m_\varphi(\cX,\cX')\subseteq \Omega^m\times {\Omega'}\phmaa{m}$ consist  of all pairs $(x,x')$ such that Duplicator has a winning strategy in the pebble game  $\sC_{m+1}(\varphi)$ on $\cX$ and $\cX'$ with initial configuration $(x,x')$. Then
\nmrt
\tm{1} $T^m_{\varphi^{}}(\cX,\cX')^*=T^m_{\varphi^{-1}}(\cX',\cX)$,
\tm{2} $[x]=[y]$ if and only if $(x,y)\in T^m_{\id_\cX}(\cX,\cX)$.
\enmrt
\elmm
\prf  Statement (1) follows from the symmetry  of $\cX$ and $\cX'$ in the game $\sC_{m+1}(\varphi)$, whereas statement~(2) is the equivalence $1\Leftrightarrow 3$ in \cite[Theorem~5.3]{CaiFI1992} for  $\cX=\cX'$ and $\varphi=\id_\cX$.
\eprf

A relationship between the pebbling game $\sC_{m+1}(\varphi)$ and the $m$-dim $\WL$-algorithm is described by the following statement proved in \cite[Lemma~3.2]{Chen2023a}.  

\lmml{131222a}
Let $\cX$ and $\cX'$ be  two ordinary coherent configurations, and let $\hat\varphi:\WL_m(\cX)\to\WL_m(\cX')$ be a bijection such that $\varphi:=\hat\varphi_2$ takes~$S$ to~$S'$.  Then
\nmrt
\tm{1} if $\hat\varphi$ is an algebraic isomorphism, then $\hat\varphi([x])=[x']$ for all $x\in\Omega^m$ and all $x'\in{\Omega'}\phmaa{m}$ such that $(x,x')\in T^m_\varphi(\cX,\cX')$,
\tm{2} if $[x]\times\hat\varphi([x])\subseteq T^m_\varphi(\cX,\cX')$ for all $x\in\Omega^m$, then $\hat\varphi$ is an algebraic isomorphism.
\enmrt
\elmm

\thrml{250124a}
Coherent configurations $\cX$ and $\cX'$ are $\WL_m$-equivalent with respect to  an algebraic isomorphism $\varphi\in\iso_{alg}(\cX,\cX')$ if and only if the binary relation $T^m_\varphi(\cX,\cX')\subseteq \Omega^m\times {\Omega'}\phmaa{m}$ has the full support.  
\ethrm
\prf
The necessity immediately follows from statement~(1) of Lemma~\ref{131222a}. To prove the sufficiency, we need the following lemma in which we put $ T^m_\varphi= T^m_\varphi(\cX,\cX')$.

\lmml{250124x}
Let $(x,x')\in T^m_\varphi$ and $y\in\Omega^m$. Then 
$$
(y,x')\in T^m_\varphi\quad \Leftrightarrow\quad [x]=[y].
$$
\elmm
\prf
First, assume that $[x]=[y]$. By Lemma~\ref{120324a}(2), we have  $(x,y)\in T^m_{\id}$, where $\id=\id_\cX$. To prove that $(y,x')\in T^m_\varphi$, we need to verify that Duplicator has a winning strategy in the game $\sC_\varphi(y,x')$ by which we denote the pebble game  $\sC_{m+1}(\varphi)$ on $\cX$ and~$\cX'$ with initial configurations $(y,x')$.  Let Spoiler chooses a set~$A'$ of points, say  in~$\cX'$. Denote by $B$ the response of Duplicator to $A'$ in the game $\sC_\varphi(x,x')$  and  by $A$  the response of Duplicator to $B$ in the game $\sC_{\id}(y,x)$. Now Duplicator responds with $A$. Obviously, $|A|=|B|=|A'|$. Spoiler places one of the pebbles  on a point $a\in A$. 

Denote by $b\in B$ the response of Duplicator to $a$ in the game $\sC_{\id}(y,x)$, and denote by $a'\in A'$ the response of Duplicator to $b$ in the game $\sC_\varphi(x,x')$. Now Duplicator responds  by placing the copy of the pebble on~$a'\in A'$. Let $\Delta=\Omega(y)\cup \{a\}$, $\Gamma=\Omega(x)\cup\{b\}$, and $\Delta'=\Omega(x')\cup\{a'\}$. Then in accordance with our definitions, the bijections
$$
f:\Delta\to \Gamma,\ y_i\mapsto x_i, \qaq
f':\Gamma\to\Delta',\ x^{}_i\mapsto x'_i, 
$$
are such that $(s_\Delta)^f=s_\Delta$ and $(s_{\Gamma^{}})^{f'}=\varphi(s)_{\Delta'}$ for all $s\in S$. It follows that the composition $f'\circ f:\Delta\to \Delta'$ takes $s_{\Delta^{}}$ to $\varphi(s)_{\Delta'}$ for all $s\in S$.  Thus  Duplicator wins and hence $(y,x')\in T^m_\varphi$.

Now assume that $(y,x')\in T^m_\varphi$. Using similar arguments for the pebble games with initial configurations $(x,x')$ and $(y,x')$, we arrive to the bijections $f:\Gamma^{}\to \Delta'$ and $f':\Delta'\to \Delta^{}$ such that $(s_{\Gamma^{}})^f=\varphi(s_{\Delta'})$ and $(s_{\Delta'})^{f'}=\varphi^{-1}(s_{\Delta^{}})$ for all $s\in S$.   The composition $f'\circ f:\Gamma\to \Delta$ takes $s_\Gamma$ to $s_\Delta$ for all $s\in S$. Again, Duplicator wins and  $(x,y)\in T^m_{\id_\cX}$, which means that $[x]=[y]$.
\eprf

Let us prove the sufficiency. From Lemma~\ref{250124x}, it follows that there is a well-defined mapping $\hat\varphi:\WL_m(\cX)\to\WL_m(\cX')$ that takes a class $[x]\in\WL_m(\cX')$ to the class $[x']$ such that $(x,x')\in T_\varphi^m$. Since the relation $T^m_\varphi$ has the full support, the mapping $\wh\varphi$ is a bijection. Finally, by the definition of the pebble game  $\sC_{m+1}(\varphi)$, the projection mapping $\hat\varphi_2$ takes any $s\in S$ to $\varphi(s)\in S'$.  By  Lemma~\ref{131222a}(2), this shows that $\wh\varphi$ is an algebraic isomorphism. Thus the coherent configurations $\cX$ and~$\cX'$ are $m$-equivalent with respect to $\varphi$. 
\eprf

\section{Circulant schemes}\label{070324c}
There is a one-to-one correspondence between circulant schemes (defined below) and the Schur rings over cyclic groups (see, e.g., \cite[Theorem~2.4.16]{CP2019}). One should not be embarrassed that most results cited in the present paper were formulated and proved in the language of Schur rings. In the rest of this section, we will change our terminology slightly to match that used for circulant schemes (and Schur rings). In this section, $G$ is a cyclic group,  $\uni_G$ denotes the identity element of~$G$, and $G_{right}\le \sym(G)$ is a group induced by right multiplications of~$G$. All unexplained details can be found in~\cite[Section~2.4]{CP2019}.

\subsection{Definitions}\label{260923w}
A coherent configuration~$\cX$ on the elements of $G$ is said to be {\it circulant scheme} over~$G$ if the group $\aut(\cX)$ contains the subgroup $G_{right}$, or equivalently 
$$
\cX\le \inv(G_{right}).
$$ 
Every relation of $\cX$ defines a Cayley graph of $G$, i.e., a circulant graph. Vice versa, if $X$ is a circulant graph, then $\WL(X)$ is a circulant scheme on the same group.

Every subgroup $H\le G$ defines an equivalence relation $e:=e(H)$ on the elements of~$G$,  the classes of which are the cosets of $H$ in~$G$. In particular, $e(\uni_G)=1_G$ and $e(G)=\bone_G$. Furthermore,  $H=\alpha e$ with $\alpha=\uni_G$, and $H\le H'$ if and only if $e(H)\subseteq e(H')$. A subgroup $H\le G$ is called an {\it $\cX$-group} if the equivalence relation $e(H)$ is a parabolic of~$\cX$. The mapping $H\mapsto e(H)$ defines  a one-to-one correspondence between the $\cX$-groups and parabolics of~$\cX$. If  $s$ is a  relation of $\cX$ and $\alpha=\uni_G$, then the subgroups $\grp{\alpha s}$ and $\{h\in G:\ (\alpha s)h=\alpha s\}$ correspond to the parabolics~$\grp{s}$ and $\rad(s)$, respectively.  Put 
$$
\secc_1(\cX)=\{\fS\in\secc(\cX):\ \uni_\fS\ne\varnothing\}.
$$ 

Let $\fS\in\secc_1(\cX)$. We associate with $\fS$ the groups $U=U(\fS)$ and $L=L(\fS)$, for which  $\fS=U/e(L)$. These are $\cX$-groups and $\fS$ can naturally be identified with cyclic group~$U/L$.
Thus $\cX_\fS$ is a circulant scheme over~$\fS$;  we say that the section~$\fS$ is \emph{trivial}  if the scheme $\cX_\fS$ is trivial.  Note that~$\fS$ is a subsection of a section $\fT\in \secc_1(\cX)$ if and only if $L(\fT)\le L$ and $U\le U(\fT)$. 

Let $U$ and $L$ be $\cX$-groups such that $L\le U$. The scheme $\cX$ satisfies the \emph{$U/L$-condition} if for any $s\in S(\cX)$ such that $s\cap e(U)=\varnothing$, we have $e(L)\le\rad(s)$. In particular, either $s\subseteq e(U)$ or $e(L)\cdot s=s\cdot e(L)=s$.

\subsection{Isomorphisms}
Let $\cX$ and $\cX'$ be circulant schemes over groups $G$ and $G'$, respectively. Under a Cayley isomorphism from $\cX$ to $\cX'$, we mean a group isomorphism $f:G\to G'$ that is a coherent configuration isomorphism. The set of all Cayley isomorphisms~$f$ is denoted by $\iso_{cay}(\cX,\cX')$, and by $\iso_{cay}(\cX)$ if $\cX=\cX'$. It was proved by Schur (see, e.g., \cite[Theorem~2.4.10]{CP2019}) that 
\qtnl{150923a}
\iso_{cay}(\cX)=\aut(G)
\eqtn 
for any circulant scheme~$\cX$ over $G$. An important step in understanding circulant schemes was done by M.~Muzychuk who proved the following theorem.

\thrml{160823c}{\rm \cite[Theorem~1.1]{Muz1994}}
Any two algebraically isomorphic circulant schemes are (combinatorially) isomorphic.
\ethrm

A complete characterization of circulant schemes was obtained by K.~Y.~Leung and K.~S.~Man, see \cite[Section~4]{Muzychuk2009}. In what follows, we will mainly use some consequences of this characterization. 

Let $\cX$ be a circulant scheme over the group $G$. From formula~\eqref{150923a}, it follows that the parabolic $\rad(s)$ does not depend on the choice of the relation~$s=r_\cX(\uni_G,g)$, where $g$ is a generator of~$G$. This parabolic is of the form $e(H)$ for a certain subgroup $H\le G$ called the \emph{radical} of~$\cX$. A section $\fS\in\secc_1(\cX)$ is said to be \emph{principal} if the radical of the scheme $\cX_\fS$ is equal to the identity subgroup~$\uni_\fS$.

\subsection{Projective equivalence}
The concept of projective equivalence comes from the lattice theory. In the case of circulant schemes, we are interested in the lattice of parabolics. It is distributive and the sections of $\secc_1(\cX)$ can be treated as  intervals in it.  A section $\fS\in\secc_1(\cX)$ is called a {\it multiple} of a section $\fT\in\secc_1(\cX)$ if
\qtnl{310114b}
L(\fT)=U(\fT)\cap L(\fS)\qaq U(\fS)=U(\fT)L(\fS).
\eqtn
Obviously, $|U(\fT)/L(\fT)|=|U(\fS)/L(\fS)|$, and $U(\fS)/U(\fT)$ is a multiple of $L(\fS)/L(\fT)$. A part of the lattice of the $\cX$-groups related to this concept is given in Fig.~\ref{fig1}.

\begin{figure}[t]
$\xymatrix@R=4pt@C=15pt@M=0pt@L=2pt{
         & \VRTB{U(\fS)} &  \\
\VRTB{L(\fS)}\ar@{-}[ur] & & \VRTB{U(\fT)}\ar@{-}[lu]     \\
		&\VRTB{L(\fT)}\ar@{-}[ur] \ar@{-}[ul] & \\
}$
\caption{$\fS$ is a multiple of $\fT$: relevant $\cX$-subgroups.}\label{fig1}
\end{figure}

The {\it projective equivalence} relation on the set $\secc_1(\cX)$ is defined to be the symmetric transitive closure of the relation ``to be multiple''; we write $\fS\sim\fT$ if the sections $\fS$ and $\fT$ are projectively equivalent, i.e., belong to the same class of projective equivalence.  

If $\fS$ is a multiple of~$\fT$, then the bijection 
\qtnl{310114b'}
f_{\fT,\fS}:\fT\to\fS,\ L(\fT)\fs\mapsto L(\fS)\fs, \ \fs \in U(\fT),
\eqtn
is a group isomorphism; hence any two projectively equivalent sections are obviously isomorphic as groups. A more general statement cited below was proved in \cite[Theorem~4.2]{Evdokimov2016}.

\lmml{261010a}
If $\fS, \fT\in\secc_1(\cX)$ and  $\fS\sim\fT$, then $f_{\fT,\fS}\in \iso_{cay}(\cX_\fT,\cX_{\fS})$.
\elmm

Let $C$ be a class of projective equivalent sections. From Lemma~\ref{261010a}, it follows that all sections of $C$ have the same order, and  if at least one section of  $C$ is trivial, then so is any section of $C$. This enables us to speak of the order of the class~$C$, and call $C$ \emph{trivial class} if it contains at least one trivial section. Put
$$
\secc_0(\cX)=\{\fS\in\secc_1(\cX):\ \fS\sim\fT\text{ for some principal }\fT\}.
$$

\subsection{Normal and quasinormal circulant schemes}\label{230923a}
The circulant scheme $\cX$ is said to be {\it normal} if the group   $G_{right}$ is normal in $\aut(\cX)$. This condition is always satisfied if the group~$G$ is of prime order unless $\cX$ is trivial and $|G|\ge 4$. A complete characterization of normal circulant schemes was obtained in \cite[Theorem~6.1]{EvdP2003}. One of the consequences of it is that in such a scheme $\cX$, every subgroup of $G$ is an $\cX$-group. Another consequence of that characterization is given in the following statement.

\lmml{220823a}
Let $\cX$ be a normal circulant scheme over a group $G=\grp{g}$. Assume that the radical of~$\cX$ is trivial. Then
\nmrt
\tm{1} the coherent configuration $\cX_{\alpha,g}$ is discrete, where $\alpha=\uni_G$,
\tm{2} for every algebraic automorphism $\varphi$ of $\cX$, we have $\iso(\cX_\alpha,\varphi)\subseteq \aut(G)$, where $\iso(\cX_\alpha,\varphi)=\{f\in\iso(\cX_\alpha):\ s^f=\varphi(s)$ for all $s\in S(\cX)\}$.
\enmrt
\elmm
\prf
Statement (1) follows from Theorem~6.1 and Corollary~6.2 of paper \cite{EvdP2003}, whereas statement~(2) follows from  Lemma~6.7(2) {\it ibid.}
\eprf

A section $\fS\in\secc_1(\cX)$ is said to be \emph{normal} if the circulant scheme $\cX_\fS$ is normal.  From Lemma~\ref{261010a}, it immediately follows that any section projectively equivalent to a normal section is normal. A subsection of a normal section is not necessarily normal, but it is always \emph{quasinormal}, i.e., is projectively equivalent to a subsection of a normal section.  A circulant  scheme  is said to be {\it quasinormal} if every its trivial section belonging to $\secc_1(\cX)$ is quasinormal.\footnote{This definition is equivalent to the original one \cite{EvdP2004a}, by the above mentioned  characterization of normal circulant schemes and  \cite[Theorem~5.1]{Evdokimov2016}.}  
An important property of a quasinormal circulant scheme $\cX$ is that the sections in $\secc_1(\cX)$ are controlled by quasinormal sections in the following sense.

\lmml{310823a1}
Let $\cX$ be a quasinormal circulant scheme and  $\fS\in\secc_1(\cX)$. Then
\qtnl{180923a}
\fS=\fS_1\times\cdots\times\fS_k\qaq \cX_\fS=\cX_{\fS_1}\otimes\cdots\otimes\cX_{\fS_k},
\eqtn
where the section  $\fS_i\in\secc_1(\cX)$ is projectively equivalent to a subsection of a principle normal section, $i=1,\ldots,k$. 
\elmm
\prf
First assume that the section $\fS$ is principal. Then the scheme $\cX_\fS$ has trivial radical. By \cite[Corollary~6.4]{EvdP2003}, there are trivial sections $\fS_1,\ldots,\fS_{k-1}\in\secc_1(\cX)$ and a trivial or principal normal section $\fS_k\in\secc_1(\cX)$ such that both equalities in~\eqref{180923a} hold true. It remains to verify that every trivial section $\fS_i$, where $i\in\{1,\ldots,k\}$, is  projectively equivalent to a subsection of a principle normal section. 

By the definition of quasinormal scheme, $\fS_i$ is projectively equivalent to a subsection of some normal section $\fT$. If the latter is principal, then we are done. Assume that the section $\fT$ is not principal. Then by \cite[Theorem~6.1 ]{EvdP2003}, the  radical of the scheme $\cX_\fT$ has order~$2$. By Corollary~6.5 {\it ibid},  there are $\cX_\fT$-groups $L$ and $U$ such that $\fT_1=\fT/L$ and $\fT_2=U$ are principal normal sections in $\secc_1(\cX_\fT)$, and any trivial $\cX_\fT$-section is a subsection of either $\fT_1$ or $\fT_2$. It follows that $\fT_1$ and $\fT_2$ are principal normal sections in $\secc_1(\cX)$, and the section~$\fS_i$ is projectively equivalent to a subsection of either $\fT_1$ or $\fT_2$.

Now let the section $\fS$ be projectively equivalent to a subsection $\fT$ of a principal section~$\fS'$. By Lemma~\ref{261010a}, the schemes $\cX_\fS$ and $\cX_\fT$ are Cayley isomorphic. Therefore without loss of generality we may assume that $\fS=\fT$, i.e., $\fS$ is a subsection of the section~$\fS'$. By the first part of the proof, the section $\fS'$ admits  decomposition~\eqref{180923a}. The restriction of this decomposition to the subsection $\fS'$ yields the required decomposition for~$\fS$.
\eprf

\section{Reduction to quasinormal schemes}\label{100324a}
\subsection{Singular classes and extensions}\label{280124a}
Let us recall some relevant facts  from papers~\cite{EvdP2004a,EvdP2012a,Evdokimov2015,Evdokimov2016}. Let $\cX$ be a circulant scheme over  a (cyclic) group $G$ and $C$ a class  of projectively equivalent sections belonging to $\secc_1(\cX)$. We say that $C$ is {\it singular} if $C$ is trivial, the order of $C$ is greater than~$2$, and $C$ contains two
sections $L_1/L_0$ and~$U_1/U_0$ such that the second is a multiple of the first and
\nmrt
\tm{S1} $\cX$ satisfies the $U_0/L_0$- and $U_1/L_1$-condition,
\tm{S2} $\cX_{U_1/L_0}=\cX_{L_1/L_0}\otimes\cX_{U_0/L_0}$.
\enmrt
By \cite[Lemma~6.2]{EvdP2012a}  the sections $L_1/L_0$ and $U_1/U_0$ are necessarily  the smallest and largest $\cX$-sections of the class~$C$, respectively. Moreover, from Theorem~4.6 {\it ibid.}, it follows that any trivial class of projectively equivalent sections of   composite order is singular. 

\lmml{240124a}
A circulant scheme $\cX$ is quasinormal if and only if each singular class of $\cX$ is of order~$3$.
\elmm
\prf
Follows from \cite[Proposition~5.3 and Theorem 5.1]{EvdP2004a}.
\eprf

Let $\fS\in\secc_1(\cX)$. The {\it $\fS$-extension} of the scheme~$\cX$ is defined to be the  smallest circulant scheme $\cX^\star\ge\cX$ such that $(\cX^\star)_\fS=\inv(\fS_{right})$, or equivalently, the scheme $(\cX^\star)_\fS$ is regular.\footnote{The $\fS$-extension always exists, because $\inv(G_{right})_\fS=\inv(\fS_{right})$.} From \cite[Theorem~4.2]{Evdokimov2016},  it follows that $\cX^\star$ does not depend on the choice of $\fS$ in the class  of $\cX$-sections projectively equivalent to~$\fS$. When $\fS$ belongs to a singular class, we say that $\cX^\star$ is a \emph{singular $\fS$-extension} or just a \emph{singular extension}.

\lmml{101213a}
Let $\cX$ be a circulant scheme, $\fS=L_1/L_0$ the smallest section in a singular class of~$\cX$, and  $\cX^\star$ the $\fS$-extension of $\cX$. Then $\rk(\cX^\star)>\rk(\cX)$ and 
\nmrt
\tm{1} the number of singular classes of $\cX^\star$ is one less than that in $\cX$, 
\tm{2} $\cX^\star$ satisfies the conditions {\rm (S1)} and {\rm (S2)} for $\cX=\cX^\star$,
\tm{3} $\{s\in S:\ s\cap e(L_1)=\varnothing\}=\{s'\in S^\star:\ s'\cap e(L_1)=\varnothing\}$,
\tm{4} $\{s\in S:\ s\subseteq e(U_0)\}=\{s'\in S^\star:\ s'\subseteq e(U_0)\}$.
\enmrt
\elmm
\prf
Statement (1) follows from \cite[Theorem~7.1(2)]{EvdP2012a}; the other three statements follow from Lemma~13.1 in~\cite{Evdokimov2016} and statements (E1), (E2) inside it.
\eprf

The theorem below immediately follows from \cite[Theorem~4.4]{Evdokimov2015}. It  is a key statement for extending algebraic isomorphisms to singular $\fS$-extensions.

\thrml{270514a}
Let $\cX$ be a circulant scheme and  $\cX^\star$ a singular $\fS$-extension of~$\cX$.  Then given $\varphi\in\aut_{alg}(\cX)$ and $\psi\in\aut_{alg}((\cX^\star)_\fS)$, there exists a unique algebraic automorphism $\varphi^\star\in\aut_{alg}(\cX^\star)$ that extends $\varphi$ and also $(\varphi^\star)_{\fS}=\psi$. 
\ethrm

\subsection{Reduction theorem} The main step in our reduction of general circulant schemes to quasinormal ones is provided by the theorem below. For $m=2$, it immediately follows from Theorem~\ref{270514a}.

\thrml{220124b}
Let  $\cX$ be a circulant scheme,  $\cX^\star$ a singular extension of~$\cX$, and $\varphi\in\aut_{alg}(\cX)$.  Assume that $\cX$ is $\WL_m$-equivalent  to itself  with respect to~$\varphi$, $m\ge 2$. Then there exists  $\varphi^\star\in\aut_{alg}(\cX^\star)$ that extends $\varphi$ and also $\cX^\star$ is $\WL_m$-equivalent to itself  with respect to~$\varphi^\star$. 
\ethrm
\prf
We keep the notation of Subsection~\ref{280124a}, concerning the sections $\fS=L_1/L_0$ and $U_1/U_0$ of a singular class of projectively equivalent $\cX$-sections, containing $\fS$. Put 
$$
e_0=e(L_0),\quad e_1=e(L_1),\quad e_2=e(U_0),\quad e_3=e(U_1).
$$
From condition (S2), it follows that if $(\alpha_0,\alpha)\in e_3$, then there exist points $\alpha_1$ and~$\alpha_2$ such that $(\alpha,\alpha_1)\in e_2$, $(\alpha,\alpha_2)\in e_1$, and $r(\alpha_0,\alpha)=r(\alpha_0,\alpha_1)\cdot r(\alpha_0,\alpha_2)$, where $r=r_\cX$. Moreover, if $(\alpha,\alpha_0)\not\in e_2$, then $\alpha_2$ can be replaced here by arbitrary point of~$\alpha_2e_0$ (because $\cX$ satisfies the $U_0/L_0$-condition, see~(S1)).

Let $G$ be the cyclic group underlying the scheme~$\cX$.
In what follows, we fix a tuple $x\in G^m$ and denote by $\varphi^\star$ the algebraic automorphism from Theorem~\ref{270514a}, that extends $\varphi$ and satisfies the condition $(\varphi^\star)_{\fS}=\id_{\cX_\fS}$. We need several auxiliary lemmas.

\lmml{020224a}
There exists a bijection $f:G\to G$ such that 
\qtnl{030224a5}
\varphi(r(x_i,\alpha))=r(x_i^f,\alpha^f),\qquad  \alpha\in G,\  i\in I,
\eqtn
where $I=\{1,\ldots,m\}$.
\elmm
\prf
Since the scheme  $\cX$ is $\WL_m$-equivalent  to itself  with respect to~$\varphi$,  Theorem~\ref{250124a} implies that there is a tuple $x'\in G^m$ such that $(x,x')$ belongs to the set $T^m_\varphi(\cX):=T^m_\varphi(\cX,\cX)$. In other words,  Duplicator has a winning strategy in the pebble game  $\sC_{m+1}(\varphi)$ on $\cX$ and $\cX'=\cX$ with initial configuration $(x,x')$. Following this strategy, for each element $\alpha\in G$ selected by the Spoiler, we associate a specific element $\alpha'\in G$ with which the Duplicator responds. In accordance with \cite[Fact~3.4.15 ]{Grohe2017}, the strategy can be defined so that the mapping $f:G\to G,$ $\alpha\mapsto\alpha'$, is a bijection. Since Duplicator wins after one round,  it follows that $x_i^f=x'_i$, $i\in I$, and equalities~\eqref{030224a5} hold.
\eprf

\lmml{020224c}
Let $f:G\to G$ be a bijection satisfying condition~\eqref{030224a5}. Assume that $\alpha\in G$ and $i\in I$ are such that $(x_i,\alpha)\not\in e_3$ or $(x_i,\alpha)\in e_2$. Then
\qtnl{040224a}
r^\star(x_i^f,\alpha_{}^f)=\varphi^\star(r^\star(x_i,\alpha)),
\eqtn
where $r^\star=r_{\cX^\star}$.
\elmm
\prf
The condition \eqref{030224a5} implies that if $(x_i,\alpha)\not\in e_3$, then $(x_i^f,\alpha^f)\not\in e_3$, and if $(x_i,\alpha)\in e_2$, then $(x_i^f,\alpha^f)\in e_2$.  By statements~(3) and~(4) of Lemma~\ref{101213a}, this implies that $r^\star(x_i,\alpha)=r(x_i,\alpha)$ and $r^\star(x_i^f,\alpha^f)=r(x_i^f,\alpha^f)$. Using condition~\eqref{030224a5} again and the fact that $\varphi^\star$ extends~$\varphi$, we obtain
$$
r^\star(x_i^f,\alpha^f)=r(x_i^f,\alpha^f)=\varphi(r(x_i,\alpha))=\varphi^\star(r(x_i,\alpha)), 
$$
as required.
\eprf

For any bijection $f:G\to G$ and a set  $\Delta\subseteq G$, denote by $M_\Delta(x,f)$ the set of all  $i\in I$ such that condition~\eqref{040224a} is satisfied for  all $\alpha\in \Delta$;  we abbreviate $M(x,f)=M_G(x,f)$ and put $m(x,f)=|M(x,f)|$. Obviously, $0\le m(x,f)\le m$. Put 
$$
m(x):=\max\{m(x,f):\ f\ \text{satisfies  condition~\eqref{030224a5}}\}.
$$ 

\lmml{020224b}
$m(x)=m$.
\elmm
\prf
By Lemma \ref{020224a}, there is a bijection $f$ satisfying condition~\eqref{030224a5}. Without loss of generality, we may assume that $m(x,f)=m(x)$. Suppose on the contrary that $m(x,f)<m$. Then there exists  an index  $i_0\in I$ for which condition~\eqref{040224a} is violated for at least one $\alpha\in G$. Let $\Delta\in G/e_3$ be such that $\alpha\in\Delta$. Then $x_{i_0}\in\Delta$ by Lemma~\ref{020224c}; in particular, $\Delta\cap\Omega(x)\ne\varnothing$. Furthermore,  by condition~\eqref{030224a5}, we have $r(x_{i_0}^f,\alpha_{\phantom{i_0}}^f)\subseteq e_3$ and hence $x_{i_0}^f$ and $\alpha^f$ belong to same class $\Delta'\in G/e_3$. Using condition~\eqref{030224a5} again, we conclude that $\Delta'=\Delta^f$ and, moreover,
$$
(\Delta\cap\Omega(x))^f=\Delta'\cap\Omega(x'),
$$
where $x'=(x_1^f,\ldots,x_m^f)$. 

Put $\fT=\Delta/e_2$ and $\fT'=\Delta'/e_2$. By the choice of the parabolics $e_3$ and $e_2$, the schemes $(\cX^\star)_{\fT^{}}$ and $(\cX^\star)_{\fT'}$ are regular. Consequently (see Subsection~\ref{110622w}), there exists an isomorphism
\qtnl{220324a}
g\in\iso((\cX^\star)_{\fT^{}}, (\cX^\star)_{\fT'},(\varphi^\star)_{\fT,\fT'}).
\eqtn
For each $i\in I$ such that $x_i\in \Delta$, we put  $\bar x_i=(x_i)_\fT$. Since the scheme $\cX^\star$ satisfies conditions~(S1) and~(S2), we make use of \cite[Theorem~6.9(1)]{EvdP2012a} to find an automorphism $h\in\aut(\cX^\star)\le\aut(\cX)$ such that 
\qtnl{220324h}
h^{G\setminus \Delta}=\id_{G\setminus\Delta}\qaq (\bar x'_i)^{\bar h}=(\bar x_i)^g,
\eqtn
where $\bar x'_i=(x'_i)_{\fT'}$ and $\bar h=h^{\fT'}$. Let $f^\star:G\to G$ be the composition of $f$ and $h$. It takes $\alpha$ to $\alpha^{fh}$ and hence
$$
\varphi(r(x_i,\alpha))=r(x_i^f,\alpha^f)=r(x_i^f,\alpha^f)^h=r(x^{fh}_i,\alpha^{fh})=r(x_i^{f^\star},\alpha^{f^\star}).
$$
Thus the bijection $f^\star$ satisfies condition~\eqref{030224a5}. Furthermore, let $i\in M_{G\setminus \Delta}(x,f)$. Then condition~\eqref{040224a} is satisfied for all $\alpha$. Indeed, if $\alpha\not\in \Delta$, then this follows because $(f^\star)^{G\setminus\Delta}=f^{G\setminus\Delta}$, whereas  if  $\alpha\in\Delta$, then this follows  from  Lemma~\ref{020224c}. Thus $i\in M(x,f^\star)$ and hence
\qtnl{050224a}
M_{G\setminus \Delta}(x,f^\star)=M_{G\setminus \Delta}(x,f).
\eqtn
To get a final contradiction, it suffices to verify that $M_\Delta(x,f^\star)=\{i\in I: x_i\in \Delta\}$. Indeed, then, in view of~\eqref{050224a}, we have 
$$
m(x)=m(x,f)=|M_{G\setminus \Delta}(x,f)|+|M_\Delta(x,f)|<
$$
$$
|M_{G\setminus \Delta}(x,f^\star)|+|M_\Delta(x,f^\star)|-1<m(x,f^\star)
$$
which contradicts the definition of $m(x)$.

It remains to prove that $M_\Delta(x,f^\star)$ contains each $i\in I$ such that $x_i\in\Delta$. Because the scheme $\cX$ is circulant, we may assume that $\Delta=U_1$. Let $\alpha\in G$. By Lemma~\ref{020224c}, we may also assume that $\alpha\in\Delta$. Then $r^\star(x_i,\alpha)\subseteq e_3$. According to remarks in the beginning of the proof for $\alpha_0=\uni_G\in U_1$, there are points $\alpha_1,\alpha_2\in\Delta$ such that $r^\star(x_i,\alpha)=r^\star(x_i,\alpha_1)\cdot r^\star(x_i,\alpha_2)$ and also $(\alpha,\alpha_1)\in e_1$ and $(\alpha,\alpha_2)\in e_2$. We have
\qtnl{220324b}
\varphi^\star(r^\star(x_i,\alpha))=\varphi^\star(r^\star(x_i,\alpha_1)\cdot r^\star(x_i,\alpha_2))=\varphi^\star(r^\star(x_i,\alpha_1))\cdot \varphi^\star(r^\star(x_i,\alpha_2)).
\eqtn
In view of condition (S1), $r^\star(x_i,\alpha_1)=\pi^{-1}(\bar r^\star(\bar x_i,\bar\alpha_1))$, where $\pi:G\to G/e_0$  is a natural surjection and  $\bar r^\star=(r^\star)_{\cX_{}}$. Furthermore, 
$r^\star(x_i,\alpha_2)=r(x_i,\alpha_2)$ by Lemma~\ref{020224c}. Thus we can continue the last equality as follows:
\qtnl{220324c}
\varphi^\star(r^\star(x_i,\alpha_1))\cdot \varphi^\star(r^\star(x_i,\alpha_2))=\varphi^\star(\pi^{-1}(\bar r^\star(\bar x_i,\bar\alpha_1)))\cdot \varphi^\star(r(x_i,\alpha_2)).
\eqtn
By definition of the isomorphism $g$, it follows that  $\varphi^\star(\pi^{-1}(\bar s))=\pi^{-1}(\bar s^g)$ for all $\bar s\in S((\cX^\star)_\fT)$.  Moreover, $\varphi$ and $\varphi^\star $ coincide on~$r(x_i,\alpha_2)$. Therefore,
\qtnl{220324d}
\varphi^\star(\pi^{-1}(\bar r^\star(\bar x_i,\bar\alpha_1)))\cdot \varphi^\star(r(x_i,\alpha_2))=\pi^{-1}(\bar r^\star(\bar x_i^g,\bar\alpha_1^g))\cdot \varphi(r(x_i,\alpha_2)).
\eqtn
Next,  using formula~\eqref{220324h} and the fact that $h$ is an automorphism of $\cX^*$, we see that $\bar r^\star(\bar x_i^g,\bar\alpha_1^g)=\bar r^\star((\bar x_i')^{\bar h},(\bar\alpha'_1)^{\bar h})=\bar r^\star(\bar x'_i,\bar\alpha'_1)$. Now,
$$
\pi^{-1}(\bar r^\star(\bar x_i,\bar\alpha_1))=\pi^{-1}(\bar r^\star(\bar x'_i,\bar\alpha'_1))=r^\star(\bar x'_i,\bar\alpha'_1)=r^\star(x_i^{f^*},\alpha_1^{f^*}),
$$
because $e_0\subseteq\rad(r^\star(\bar x'_i,\bar\alpha'_1))$. Finally, taking into account that     $\varphi(r(x_i,\alpha_2))=r^\star(x_i^{f^\star},\alpha^{f^\star})$ by Lemma~\ref{020224c}, we can continue equality \eqref{220324d} as follows:
\qtnl{220324e}
\pi^{-1}(\bar r^\star(\bar x_i^g,\bar\alpha_1^g))\cdot \varphi(r(x_i,\alpha_2))=r^\star(x_i^{f^\star},\alpha_1^{f^\star})\cdot r(x_i^{f^\star},\alpha_2^{f^\star})=r^\star(x_i^{f^\star},\alpha_{}^{f^\star}).
\eqtn
Collecting together formulas \eqref{220324b}, \eqref{220324c}, \eqref{220324d}, and~\eqref{220324e}, we obtain $\varphi^\star(r^\star(x_i,\alpha))=r^\star(x_i^{f^\star},\alpha_{}^{f^\star})$, as required.
\eprf

Let us complete the proof of Theorem~\ref{220124b}. From Lemma~\ref{020224b}, it follows that there exists a bijection $f:G\to G$ for which $m(x,f)=m$. Let us verify that  $(x,x')\in T^m_{\varphi^\star}(\cX^\star)$, where  $x'=(x_1^f,\ldots,x_m^f)$. In other words, we need to  verify that  Duplicator has a winning strategy in the pebble game  $\sC_{m+1}(\varphi^\star)$ on the scheme $\cX^\star$ with initial configuration $(x,x')$.  Let Spoiler choose a set~$A'$ of points in, say~$\cX$. Duplicator responds  with  set $A=A^{f^{-1}}$. Obviously,  $|A|=|A'|$. Now when Spoiler puts a pebble on a point $a\in A$, Duplicator puts the corresponding pebble to the point $a^f$, and wins by the choice of~$f$.

Since the $m$-tuple $x$ was arbitrary, the last paragraph shows that  under the theorem condition, we have $ \pr_1 T_{\varphi^\star}^m(\cX^\star)=G^m$. Applying the same argument to the algebraic automorphisms $\varphi^{-1}$ and $(\varphi^\star)^{-1}$, we obtain 
$$
\pr_2 T_{\varphi^\star}^m(\cX^\star)=\pr_1 T_{(\varphi^\star)^{-1}}^m(\cX^\star)=G^m.
$$
Thus the binary relation $T_{\varphi^\star}^m(\cX^\star)$ has the full support. It follows that then  $\cX^\star$ is $\WL_m$-equivalent to itself  with respect to~$\varphi^\star$ by Theorem~\ref{250124a}.
\eprf

\section{Point extensions of a circulant scheme}
Throughout this section, $\cX$ is a circulant scheme over a cyclic group $G$. For any point tuple $x$ of~$\cX$, we put $\cX_{x,\fS}:=(\cX_x)_\fS$ for all  $\fS\in\secc(\cX)$.

\subsection{Auxiliary lemmas}
According to our definitions, each  element $\fs$ of a section  $\fS\in\secc(\cX)$ can be treated in two ways: as  a class of the parabolic $e(L)$ with  $L=L(\fS)$, and as a coset of $L$ in the group $U(\fS)$. In the first case, $\fs$ is a subset of~$G$, whereas in the second one, $\fs$ is an element of $\fS$. Each time it will be clear from the context, in what capacity the element~$\fs$ is considered. The following statement is straightforward.

\lmml{050923d1}
Let a section $\fS\in\secc_1(\cX)$ be a multiple of a section $\fT\in\secc_1(\cX)$. Then for any $\fs\in\fS$ and any $\ft\in\fT$, we have
$$
\ft^{f_{\fT,\fS}}=\fs\quad\Leftrightarrow\quad \fs\supseteq \ft.
$$
\elmm

The sections of the scheme $\cX$ are also sections  in any  extension of~$\cX$. The following lemma provides simple tools to work with sections of the point extensions of $\cX$.

\lmml{090923a}
Let $\fs\in\fS\in\secc_1(\cX)$, $x$ a point tuple of~$\cX$, and $\uni_G\in\Omega(x)$. Then
\nmrt
\tm{1} $\{\fs\}\in F(\cX_{x,\fS})$ if and only if $\fs\in F(\cX_x)^\cup$,
\tm{2} if  $\fS\sim\fT$, then  $\{\fs\}\in F(\cX_{x,\fS})$ if and only if $\{\fs^{f_{\fS,\fT}}\}\in F(\cX_{x,\fT})$. 
\enmrt
\elmm
\prf
The condition $\uni_G\in\Omega(x)$ implies that every $\cX$-subgroup of $G$ is a homogeneity set of the coherent configuration~$\cX_x$. In particular, $U(\fS)\in F(\cX_x)^\cup$.  Thus statement~(1) follows from the fact that the quotient map $U(\fS)\to U(\fS)/L(\fS)$ induces an (inclusion preserving) surjection from the homogeneity sets of the coherent configuration $(\cX_x)_{U(\fS)}$ to the homogeneity sets of the coherent configuration  $\cX_{x,\fS}$.

To prove statement~(2), we assume without loss of generality that one of the sections $\fS$ and $\fT$ is a multiple of the other, say $\fT$ is a multiple of $\fS$. Let $\ft=\fs^{f_{\fS,\fT}}$. By Lemma~\ref{050923d1}, we have $\fs\subseteq\ft\subseteq U(\fT)$  and hence
$$
\ft=L(\fT)g_\ft\qaq \fs=L(\fS)g_\fs
$$
for some elements $g_\ft\in U(\fT)$ and $g_\fs\in U(\fS)$. It follows that $g_\fs$ belongs to the coset $(U(\fS)\cap L(\fT))g_\ft$ of the group $U(\fS)\cap L(\fT)=L(\fS)$.  Thus, 
\qtnl{100124a}
U(\fS)\cap\ft=L(\fS)g_\fs=\fs.
\eqtn

Assume that $\{\fs^{f_{\fS,\fT}}\}=\{\ft\}$ belongs to $F(\cX_{x,\fT})$. By statement~(1), this implies that $\ft\in  F(X_x)^\cup$. Since also $U(\fS)\in F(\cX_x)^\cup$, we conclude  by~\eqref{100124a} that  $\{\fs\}=\ft\cap U(\fS)\in F(\cX_x)^\cup$.  Using  statement~(1) again, we obtain $\{\fs\}\in F(\cX_{x,\fS})$, as required. 

To prove the converse statement, assume that $\{\fs\}\in F(\cX_{x,\fS})$. Note that $\fs=L(\fS)g_\fs\subseteq L(\fT)g_\fs=L(\fT)g_\ft$ implying that
$$
1_\fs\,\cdot\, e(L(\fT))=L(\fT)g_\ft\times L(\fT)g_\ft.
$$ 
The support of the relation on the right-hand side is $\ft$. By the assumption,  the relation~$1_\fs$ and hence the product $1_\fs\,\cdot\, e(L(\fT))$ is a relation of $\cX_x$. Thus its support~$\ft$  belongs to  the set~$F(\cX_x)^\cup$. By statement~(1), this implies that $\{\ft\}\in F(\cX_{x,\fT})$.
\eprf

\subsection{Base tuples}
A  point tuple $x$ of a circulant scheme~$\cX$   is called a {\it base} tuple for this  scheme if $\uni_G\in\Omega(x)$ and for each section $\fS\in\secc_1(\cX)$ of prime power order, there is  $g\in\Omega(x)$ such that $\fS=\grp{g_\fS}$. A trivial example of a base tuple is given by any $x$ for which $\Omega(x)=G$. 

\lmml{160923a}
For any circulant scheme of degree $n$, there exists a base $m$-tuple with $m\le \Omega(n)+1$. 
\elmm
\prf
 Denote by $\cH$  the set of all subgroups of $G$ of prime power order. For each $H\in \cH$, we choose an arbitrary generator $g_H$ of the group~$H$. Then any point tuple $x$ with  $\Omega(x)=\{g_H:\ H\in\cH\}$ is a base tuple for~$\cX$: indeed, if $\fS\in\secc_1(\cX)$  is of prime power order $p^k$ for some $k\ge 1$, then the set $U(\fS)\setminus L(\fS)$ contains $g_H$, where $H$ is the Sylow  $p$-subgroup of~$U(\fS)$. Since the group $G$ is cyclic, we have $|\cH|\le\Omega(n)+1$, as required.
 \eprf

The following statement shows that if $x$ is a base tuple for a quasinormal scheme~$\cX$, then the set $\Omega(x)$ controls, in a sense, the sections belonging to the  set $\secc_0(\cX)$.

\thrml{170823a1}
Let  $x$ be a base tuple for a quasinormal circulant scheme~$\cX$.  Then for any section  $\fS\in \secc_0(\cX)$, the  coherent configuration $\cX_{x,\fS}$ is discrete.
\ethrm
\prf
First, let $\fS\in\secc_0(\cX)$ be a normal principle section. Then for any prime~$p$ dividing $|\fS|$, the Sylow $p$-subgroup $\fS_p$ of $\fS$ and its complement subgroup $\fS_{p'}$ are $\cX_\fS$-groups. Since $x$ is a base tuple, there exists $g\in\Omega(x)$ such that $U(\fS_p)=\grp{g}$. Put $\fs_p=g_\fS$. Then 
$$
\grp{\fs_p}=\fS_p\qaq \{\fs_p\}\in F(\cX_{x,\fS}),
$$ 
where the latter holds true, because $g\in\Omega(x)$ and hence the singleton $\{g\}$ is a fiber of the coherent configuration $\cX_g \le \cX_x$. By the same reason, $\{\uni_\fS\}\in F(\cX_{x,\fS})$, and hence $\fS_{p'}$ is a homogeneity set of $\cX_{x,\fS}$. It easily follows that
\qtnl{270224a}
\fs_p\fS_{p'}\in F(\cX_{x,\fS})^\cup.
\eqtn
Note that the intersection of  the sets $\fs_p\fS_{p'}\subseteq \fS$ over all $p$ dividing~$|\fS|$ is equal to a singleton $\{\fs\}$, where $\fs=\prod_p \fs_p\in\fS$.  In view of~\eqref{270224a}, this singleton is a homogeneity set of the coherent configuration $\cX_{x,\fS}$; in particular,  $\{\fs\} \in F(\cX_{x,\fS})$. On the other hand, the element $\fs$ generates the group $\fS$ (see above). Thus, the  coherent configuration
$$
\cX_{x,\fS}\ge (\cX_\fS)_{\uni_\fS,\fs}
$$
is discrete by Lemma~\ref{220823a}(1).

Now, let $\fS\in\secc_0(\cX)$ be an arbitrary section. It is of the form described in  Lemma~\ref{310823a1}. In the notation of that lemma, from \cite[Theorem~3.2.17]{CP2019}, it follows that 
$$
\cX_{x,\fS}\ge (\cX_{x,\fS})_{\fS_1}\otimes\cdots\otimes(\cX_{x,\fS})_{\fS_k}=
\cX_{x,\fS_1}\otimes\cdots\otimes\cX_{x,\fS_k}.
$$
It suffices to verify that the coherent configuration $\cX_{x,\fS_i}$ is discrete, $i=1,\ldots,k$. Recall that the section $\fS_i$  is projectively equivalent to a subsection $\fT$ of some principal normal section $\fT'\in\secc_0(\cX)$. By above, the coherent configuration $\cX_{x,\fT'}$ and hence the coherent configuration $\cX_{x,\fT}=(\cX_{x,\fT'})_\fT$ is discrete. By Lemma~\ref{090923a}(2), this implies that all fibers of the coherent configuration $\cX_{x,\fS}$ are singletons, i.e., this coherent configuration  is discrete.
\eprf

\section{Extendable algebraic isomorphisms and multipliers}

\subsection{Extendable algebraic isomorphisms} 
Let $\cX$ and $\cX'$ be coherent configurations and $x$ a point tuple of~$\cX$.   An algebraic isomorphism $\varphi\in\iso_{alg}(\cX,\cX')$ is said to be {\it  extendable} at $x$ if $\varphi$ has the $(x,x')$-extension for some point tuple~$x'$ of~$\cX'$. Let $m$ be a nonnegative integer. If the algebraic isomorphism $\varphi$ is extendable at every $m$-tuple, then $\varphi$ is said to be  {\it $m$-extendable}. Thus the usual algebraic isomorphisms are $0$-extendable, whereas the sesquiclosed algebraic isomorphisms defined in \cite{Ponomarenko2022a} are $1$-extendable.

\lmml{240823e}
Let  $\cX$ and $\cX'$ be schemes, $\varphi\in\iso_{alg}(\cX,\cX')$, and $m$ a nonnegative integer. Assume that $\varphi$ is $m$-extendable. Then 
\nmrt
\tm{i} $\varphi$ is  $k$-extendable for all $0\le k\le m$,
\tm{ii} if  $\fS\in\secc(\cX)$,  $\fS'\in\secc(\cX')$, and  $\varphi_{\fS,\fS'}\in\iso_{alg}(\cX^{}_{\fS^{}},\cX'_{\fS'})$ is defined, then  $\varphi_{\fS,\fS'}$ is $m$-extendable.
\enmrt
\elmm
\prf 
Let $y$ be a point $k$-tuple of~$\cX$. Choose any $m$-tuple $x$ such that $\pr_k(x)=y$. By the lemma condition, there is an $(x,x')$-extension $\psi$ of the algebraic isomorphism~$\varphi$ for some point  $m$-tuple~$x'$ of $\cX'$. The restriction of  $\psi$ to the coherent configuration~$\cX_y$ is obviously the $(y,y')$-extension of~$\varphi$. This proves statement~(i). Statement~(ii) follows from Lemma~\ref{030923h}.
\eprf

The concept of $m$-extendable algebraic isomorphism can be used to establish a necessary condition for two schemes to be $\WL_m$-equivalent. More exactly, the following statement immediately follows from \cite[Theorem~3.7]{Ponomarenko2022a}.

\lmml{310723a}
Let $m\ge 2$ be an integer. Assume that two schemes are $\WL_m$-equivalent with respect to an algebraic isomorphism. Then it is $(m-2)$-extendable.
\elmm

\subsection{Multipliers of circulant schemes}
In the rest of this section we show that if~$\cX$ is a quasinormal circulant scheme, then every $m$-extendable algebraic automorphism of $\cX$ for sufficiently large~$m$ is  induced by an isomorphism. We begin with an auxiliary statement showing that in this case the algebraic automorphism in question is associated with a special system of group automorphisms parametrized by the $\cX$-sections belonging to the set $\secc_0(\cX)$.

\lmml{160823a}
Let $\cX$ be a quasinormal circulant scheme over a group~$G$, and $\varphi$ an algebraic automorphism of $\cX$. Assume that $\varphi$ is extendable  at a base tuple $x\in G^m$.  Then there exists a family
\qtnl{240823w}
\fM=\fM(\varphi,x)=\{\sigma_\fS\in \aut(\fS):\ \fS\in\secc_0(\cX)\}
\eqtn
such that  for all sections $\fS,\fT\in\secc_0(\cX)$ the following conditions are satisfied:
\nmrt
\tm{M1} $\sigma_\fS\in\iso(\cX_\fS,\varphi_\fS)$,
\tm{M2} if $\fT\succeq \fS$, then $(\sigma_\fT)^\fS=\sigma_\fS$,
\tm{M3} if $\fT\sim \fS$, then $(\sigma_{\fT})^{f_{\fT,\fS}}=\sigma_\fS$,
where $f_{\fT,\fS}$ is the mapping \eqref{310114b'}.
\enmrt
\elmm
\prf
Let $x'\in G^m$ be an arbitrary $m$-tuple such that the algebraic automorphism $\varphi$ can be extended to an algebraic isomorphism  $\psi\in\iso_{alg}(\cX_{x^{}},\cX_{x'})$. Let $\fS\in\secc_0(\cX)$.  By Theorem~\ref{170823a1}, the coherent configuration $\cX_{x,\fS}$ is discrete. Consequently, the algebraic isomorphism $\psi_\fS: \cX_{x,\fS}\to \cX_{x',\fS}$ is induced by a (uniquely defined) permutation $\sigma_\fS\in\sym(\fS)$ such that 
\qtnl{050923a}
\psi_\fS(1_\fs)=1_{\fs^\sigma}\ \,\text{for all}\ \,\fs\in\fS,
\eqtn
where $\sigma=\sigma_\fS$.
First we verify that this permutation satisfies the conditions~(M1), (M2), (M3).

To prove (M1), let $\fS\in\secc_0(\cX)$. Since $\cX_{x,\fS}\ge \cX_\fS$ and $\psi$ extends $\varphi$, we conclude that $\psi_\fS$ extends $\varphi_\fS$, and hence $\iso(\cX_{x,\fS},\psi_\fS)\subseteq \iso(\cX_\fS,\varphi_\fS)$. On the other hand,  $\sigma_\fS\in\iso(\cX_{x,\fS},\psi_\fS)$.  Thus,  $\sigma_\fS\in\iso(\cX_\fS,\varphi_\fS)$. 

To prove (M2), let $\fS,\fT\in\secc_0(\cX)$ be such that $\fT\succeq \fS$. Then, obviously, $(\psi_\fT)_\fS=\psi_\fS$. Take any $\fs\in\fS$   and choose any $\ft\in\fT$ such that $\ft_\fS=\fs$. Then, by formula~\eqref{050923a}, we have
$$
1_{(\fs)^{\mu^\fS}}=1_{{(\ft_\fS)}^{\mu^\fS}}=(1_{\ft^\mu})_\fS=
(\psi_\fT(1_\ft))_\fS=[(\psi_\fT)_\fS](1_\fs)=\psi_\fS(1_\fs)=1_{\fs^{\sigma}},
$$
where $\mu=\sigma_\fT$. Hence, $\fs^{\mu^\fS}=\fs^\sigma$ for all $\fs\in\fS$, whence $(\sigma_\fT)^\fS=\mu^\fS=\sigma=\sigma_\fS$.

To prove (M3), let $\fS,\fT\in\secc_0(\cX)$ be such that $\fT\sim \fS$. Without loss of generality, we may assume that $\fS$ is a multiple of $\fT$: indeed, the general case reduces to the case when one of the sections $\fS$ and $\fT$ is a multiple of the other, and then we make use of the fact that $f^{}_{\fS,\fT}=(f_{\fT,\fS})^{-1}$. Thus we need to verify that 
$$
\ft^{\sigma_\fT f_{\fT,\fS}}=\fs^{\sigma_\fS}\ \,\text{for all}\ \,\ft\in\fT,
$$
where $\fs=\ft^{f_{\fT,\fS}}$. By Lemma~\ref{090923a}(1), the set $\ft$ treated as a subset of~$G$ is a homogeneity set of the coherent configuration $\cX_x$. Therefore the set $\ft^\psi$ is also a homogeneity set of  $\cX_x$. Furthermore,
$$
(\ft^\psi)_\fT=(\ft_\fT)^{\psi_\fT}=\{\ft\}^{\psi_\fT}=\{\ft^{\sigma_\fT}\}.
$$
Consequently, $\ft^\psi\subseteq \ft^{\sigma_\fT}$. The converse inclusion can be verified analogously. Thus, $\ft^\psi=\ft^{\sigma_\fT}$ and, similarly, $\fs^\psi=\fs^{\sigma_\fS}$. Now, since $\fs=\ft^{f_{\fT,\fS}}$,  Lemma~\ref{050923d1} implies that $\ft\supseteq\fs$. Consequently,  $\ft^\psi \supseteq\fs^\psi$, which by above shows that  $\ft^{\sigma_\fT}\supseteq \fs^{\sigma_\fS}$. Thus, $\ft^{\sigma_\fT f_{\fT,\fS}}=\fs^{\sigma_\fS}$ again by Lemma~\ref{050923d1}.

It remains to verify that  $\sigma_\fS\in \aut(\fS) $ for all $\fS\in\secc_0(\cX)$. This follows from condition~(M1) and Lemma~\ref{220823a}(2) if  $\fS$ is a principal normal section. Let $\fS$ be a subsection of a principal  normal  section $\fS'\in \secc_0(\cX)$ such that $\fS\sim \fT\preceq \fS'$ for some section $\fT\in\secc_1(\cX)$. By condition (M2), we have 
$$
\sigma_\fT=(\sigma_{\fS'})^\fT\in\aut(\fS')^\fT=\aut(\fT). 
$$
By condition (M3) and Lemma~\ref{261010a}, this implies that
$$
\sigma_\fS=f_{\fS,\fT}\,\sigma_\fT\, {f_{\fS,\fT}}^{-1}\in f_{\fS,\fT}\,\aut(\fT)\, {f_{\fS,\fT}}^{-1}=\aut(\fS).
$$
Finally, assume that the section $\fS\in\secc_0(\cX)$ is arbitrary. It is of the form~\eqref{180923a}, see  Lemma~\ref{310823a1}. In the notation of that lemma, from the condition (M2) and the above, it follows that 
$$
\sigma_\fS=(\sigma_\fS)^{\fS_1}\cdot\ldots\cdot(\sigma_\fS)^{\fS_k}=
\sigma_{\fS_1}\cdot\ldots\cdot\sigma_{\fS_k}\in\aut(\fS_1)\times\cdots\times\aut(\fS_k)=\aut(\fS),
$$
as required.
\eprf

From the proof of Lemma~\ref{160823a}, it follows that the automorphisms $\sigma_\fS$ in the statement of that lemma depend not only on the $m$-tuple $x$ but also on  the choice of the $m$-tuple $x'$. However, this dependence is not essential for further exposition.

At this point, we make use of some results proved in papers \cite{Evdokimov2016} and~\cite{Evdokimov2015}. We apply these results to quasinormal circulant schemes.  It should be noted that any such scheme is \emph{quasidence} in terms of that papers. Indeed, every trivial $\cX$-section is projectively equivalent to a subsection of a normal section and hence is of prime order. Moreover, the quasinormality implies also that the set  $\secc_0(\cX)$ defined in our paper coincides with the set $\fS_0(\cX)$  used in that papers. Thus any family~\eqref{240823w} satisfying the conditions (M1), (M2),  and~(M3) is (in terms of that papers) an \emph{$\cX$-multiplier}, see \cite[Definition~5.1]{Evdokimov2015}. Based on the theory of $\cX$-multipliers and using Lemma~\ref{160823a}, we can prove the main result of this section.

\thrml{170823h}
Let $\cX$ be a quasinormal circulant scheme. Then any algebraic automorphism of $\cX$, that is extendable  at a base tuple for  $\cX$,  is induced by an isomorphism.
\ethrm
\prf
Assume that $\varphi\in\aut_{alg}(\cX)$ is extendable  at  a base tuple~$x$ for~$\cX$. By Lemma~\ref{160823a}, there exists an $\cX$-multiplier $\fM=\fM(\varphi,x)$. Denote by $\cX_0$ the \emph{coset closure} of the scheme~$\cX$; there is no need to cite here the exact definition that was given in~\cite{Evdokimov2016}. For our purposes, it suffices to know just that $\cX_0$ is a circulant scheme and $\cX_0\ge \cX$. In accordance with Lemma~5.2 and Theorem~3.4 proved in \cite{Evdokimov2015}, there exists an algebraic automorphism $\varphi_0\in\aut_{alg}(\cX_0)$ such that
\qtnl{240823m}
\{\sigma_\fS\}=\iso(\cX_{0,\fS},\varphi^{}_{0,\fS})\quad\text{for all}\ \,\fS\in\secc_0(\cX),
\eqtn
where $\sigma_\fS$ is the component of $\fM$ at the section $\fS$, $\cX_{0,\fS}=(\cX_0)_\fS$,  and $\varphi_{0,\fS}=(\varphi_0)_\fS$. By \cite[Theorem~1.1]{Evdokimov2015}, it suffices to verify that $\varphi_0$ extends~$\varphi$. To this end,  let $s$ be a basis relation of the scheme~$\cX$. Denote by $\fS$ the section of $\secc_0(\cX)$ such that $U(\fS)$ and $L(\fS)$ are  the subgroups of $G$ corresponding to the parabolics~$\grp{s}$ and $\rad(s)$ in the sense of Subsection~\ref{260923w}.  From formula~\eqref{240823m}, it follows that 
$$
\varphi_0(s)_\fS=\varphi_{0,\fS}(s_\fS)=(s_\fS)^{\sigma_\fS}=\varphi(s_\fS).
$$
Taking into account that $\varphi_0(\rad(s))=\rad(s)=\varphi(\rad(s))$, we conclude that $\varphi_0(s)=\varphi(s)$. 
\eprf

\section{Proof of Theorem~\ref{260723a}} 
Let $X$ be a circulant graph of order~$n$. Then the coherent configuration $\cX=\WL(X)$ is a circulant scheme. In accordance with \cite[Subsection~3.2]{Chen2023a}, the $\WL$-dimension of $X$ is equal to $\dim_{\scriptscriptstyle\WL}(\cX)$. We estimate the  last number in two steps. First, we use Theorem~\ref{220124b} to reduce  the general  circulant schemes to quasinormal schemes. 

\thrml{210124b}
For any circulant scheme $\cX$, there exists a quasinormal circulant scheme $\cX'$ (of the same degree) such that $\dim_{\scriptscriptstyle\WL}(\cX)\le \dim_{\scriptscriptstyle\WL}(\cX')$.
\ethrm
\prf
If the scheme $\cX$ is quasinormal, then it is nothing to prove. Otherwise,  $\cX$ has a singular class of order at least $4$ (Lemma~\ref{240124a}). Let $\cX^\star$ be a singular extension of~$\cX$. Since $\cX^\star$ is circulant and $\rk(\cX^\star)>\rk(\cX)$ (Lemma~\ref{101213a}), we may assume by induction on $\rk(\cX)$  that 
\qtnl{080524a}
\dim_{\scriptscriptstyle\WL}(\cX^\star)\le \dim_{\scriptscriptstyle\WL}(\cX'):=m
\eqtn
for some quasinormal circulant scheme $\cX'$ with the same degree as $\cX^\star$.  

Now  let $\cY$  be a circulant scheme  $\WL_m$-equivalent to~$\cX$ with respect to an algebraic isomorphism~$\varphi\in\iso_{alg}(\cX,\cY)$. By the Muzychuk theorem (Theorem $\ref{160823c}$), we may assume that \qtnl{270224v}
\cY=\cX\qaq \varphi\in\aut_{alg}(\cX).
\eqtn 
Then $\cX$ and $\varphi$ satisfy the condition of Theorem~\ref{220124b}. Hence there exists  an algebraic automorphism $\varphi^\star\in\aut_{alg}(\cX^\star)$ that extends $\varphi$ and such that $\cX^\star$ is $\WL_m$-equivalent to itself  with respect to~$\varphi^\star$. By virtue of~\eqref{080524a}, this shows that  
$$
\iso(\cX^\star,\varphi^\star)\ne\varnothing.
$$  
Since $\iso(\cX^\star,\varphi^\star)\subseteq\iso(\cX,\varphi)$, we conclude that $\iso(\cX,\varphi)\ne\varnothing$. This proves the required inequality $\dim_{\scriptscriptstyle\WL}(\cX)\le m=\dim_{\scriptscriptstyle\WL}(\cX')$.
\eprf

To complete the proof, let $m=\Omega(n)+3$. By Theorem~\ref{210124b}, we may assume   that the scheme $\cX$ is quasinormal. Let $\cY$ be an arbitrary circulant scheme  $\WL_m$-equivalent to~$\cX$ with respect to an algebraic isomorphism~$\varphi\in\iso_{alg}(\cX,\cY)$. As above, we may also assume that relations~\eqref{270224v} hold. Then the algebraic isomorphism~$\varphi$ is  $(m-2)$-extendable by Lemma~\ref{310723a}. From Lemma~\ref{160923a}, it follows that there exists a base $(m-2)$-tuple  for~$\cX$. Since $\varphi$ is   extendable at this tuple, Theorem~\ref{170823h} implies that $\iso(\cX,\varphi)\ne\varnothing$. Thus, $\dim_{\scriptscriptstyle\WL}(\cX)\le m=\Omega(n)+3$, as required.

\bibliographystyle{amsplain}

\end{document}